\documentclass[reqno,12pt]{amsart}

\usepackage{amsthm,amsfonts}
\usepackage[centertags]{amsmath}
\usepackage[a4paper,margin=34mm]{geometry}

\usepackage{graphics}
\usepackage{epsfig}

\usepackage{bm}

\newtheorem{theorem}{Theorem}
\newtheorem{proposition}{Proposition}

\newtheorem{definition}{Definition}

\newtheorem{remark}{Remark}

\begin{document}

\title[]{Discrete Capacity and Higher-order Differences of Two-state Markov Chains}

\author[A. Yu. Shahverdian]{A. Yu. Shahverdian}

\address{Institute for Informatics and Automation Problems of NAS RA}

\email{svrdn@yerphi.am}



\date{\today}

\maketitle

\begin{abstract}
The paper studies the time-homogeneous two-state Markov chains; the states are assumed to be binary symbols 0 and 1. The higher-order absolute differences taken from progressive states of a given chain are considered. A discrete capacity of subsets of natural series is defined and a limiting theorem for these differences, formulated in terms of Wiener criterion type relation, is presented.
\end{abstract}

\section{\textbf{Introduction}}\label{S1}

In this paper an application of the suggested in \cite{1}-\cite{7} difference analysis to studying the two-state Markov chains is presented. The difference analysis is a method for studying irregular and random time series, based on consideration of higher-order absolute differences taken from the series' progressive terms. This method allowed us to reveal some new aspects in dynamical systems: e.g., the higher-order-difference version for Lyapunov exponent  \cite{3} and bistability of higher-order differences, taken from periodic time series \cite{6}, have been established.

We study time-homogeneous Markov chains $\bm{\xi}= (\xi_{n})_{n=0}^{\infty}$, whose state space $X =\{x\}$  consists of two different items; more precisely, we suppose that $X=\{0,1\}$, that is, each component $\xi_n$ of $\bm{\xi}$ (which describes the chain at the moment $n$) is a random binary variable.

The main result of this paper, Theorem~\ref{T1}, is a limiting theorem for such chains: it assert the existence of the limit of $k$th order absolute differences taken from progressive terms of a given series $(\xi_n)_{n=0}^{\infty}$, when $k$ converges to $\infty$ remaining on "large" subsets $E\subseteq \mathbb{N}$ of natural series $\mathbb{N}$. The "size" of such sets $E$ is described in terms of some discrete capacity: such sets $E$ are  {\em thick} sets, defined by means of Wiener criterion type relation from potential theory (see, e.g., \cite {71} and \cite{8}). The limiting process, whose existence asserts Theorem 1, is the equi-distributed random sequence.

The paper consists of three sections. The next Section~\ref{S2} describes the statement of the considered problem, in Section~\ref{S3} we present the definitions of discrete capacity, {\em thin} and {\em thick} sets, and formulate our Theorem \ref{T1}.

\section{\textbf{Statement of the problem}}\label{S2}

Let us explain the statement of the problem which we study. Let
$$
\bm{\xi} = (\xi_{0}, \xi_{1}, \ldots, \xi_{n}, \ldots)
$$
be some random sequence whose components $\xi_{n}$ take binary values $x$ from $X = \{0, 1\}$
with some positive probabilities $p_{n}(x)$, $P(\xi_{n}= x) = p_{n}(x)$ ($p_{n}(0)+p_{n}(1)=1$). Then $k$th order ($k\geq 0$) absolute differences $\xi_{n}^{(k)}$, defined recurrently as: $\xi_{n}^{(0)} \equiv \xi_{n}$ and
$$
\xi_{n}^{(k)} = |\xi_{n+1}^{(k-1)} - \xi_{n}^{(k-1)}|  \qquad (n\geq 0),
$$
also take binary values with some probabilities $p_{n}^{(k)}(x)$,
$$
P(\xi_{n}^{(k)}= x) = p_{n}^{(k)}(x) \qquad  (p_{n}^{(k)}(0)+p_{n}^{(k)}(1)=1);
$$
hence, one can consider $k$th order  difference random binary sequence
$$
\bm{\xi}^{(k)} = (\xi_{0}^{(k)}, \xi_{1}^{(k)}, \ldots, \xi_{n}^{(k)}, \ldots ).
$$
We are interested in existence of the limit of $\bm{\xi}^{(k)}$ when $k$ goes to infinity. Let some infinite $\Lambda \subseteq \mathbb{N}$ be given; we say that $\bm{\xi}^{(k)}$ converge to a random binary sequence $\bm{\xi}^{(\infty)}$, if $p_{n}^{(k)}(x)$ ($n\in \mathbb{N}$, $x\in X$) tend to some numbers $p_{n}^{(\infty)}(x)$ ($p_{n}^{(\infty)}(0)+p_{n}^{(\infty)}(1)=1$) as $k\to \infty$ and $k\in \Lambda$ (convergence by probability on $\Lambda$). Given $\Lambda$ the limiting process
$$
\bm{\xi}^{(\infty)} = \bm{\xi}^{(\infty)}_{\Lambda} = (\xi_{0}^{(\infty)}, \xi_{1}^{(\infty)}, \ldots, \xi_{n}^{(\infty)}, \ldots)
$$
(so-called {\em partial limit}) is defined as random sequence, whose components $\xi_{n}^{(\infty)}$ take the values $x\in X$ with the probabilities $p_{n}^{(\infty)}(x)$.

We study time-homogeneous Markov chains $\bm{\xi}$, that is, when for $x$, $x_i$, $y \in X$
\begin{equation}\label{e1}
P(\xi_{n}=y|\xi_{n-1}=x, \xi_{n-2}=x_{1}, \ldots, \xi_{0}=x_{n-1}) = P(\xi_{n}=y|\xi_{n-1}=x)
\end{equation}
(Markov property) and there is some function $\pi(x,y)$ on $X\times X$ such that
\begin{equation}\label{e2}
P(\xi_{n}=y|\xi_{n-1}= x) = \pi(x,y) \quad \text{for} \ n\geq 1 \ \text{and} \ x,y \in X
\end{equation}
(homogeneity). Some computations testify, that if for such $\bm{\xi}$ an infinite $\Lambda \subseteq \mathbb{N}$ is chosen arbitrarily, the limiting process $\bm{\xi}^{(\infty)}_{\Lambda}$ may not exist; on the other hand, a theorem announced in \cite{7} asserts that if $\Lambda = \{2^{m}-1: m\geq 0\}$, then $\bm{\xi}^{(\infty)}_{\Lambda}$ exists. The problem which studies the present paper is the following (descriptively): how "large" can be the sets $\Lambda\subseteq \mathbb{N}$ which permit the existence of $\bm{\xi}^{(\infty)}_{\Lambda}$, and how their "size" can be described? This paper considers the chains for which
\begin{equation}\label{E}
\pi(x,y) \neq 0, \qquad \pi(0,0) \neq \pi(1,1), \quad \text{and} \quad \pi(0,0)+\pi(1,1) \neq 1.
\end{equation}

We claim that for time-homogeneous binary Markov chains the problem stated is resolved in terms of some discrete capacity defined on $2^{\mathbb{N}}$ and corresponding {\em thin} ({\em fine}) and {\em thick} sets. The capacity $\mathcal{C}$, considered here, is a modification of the discrete capacity used in \cite{4}. The solution to our problem is given by Theorem~\ref{T1}, which is formulated in terms of thick sets, defined by means of well-known in potential theory Wiener criterion type relation.

\section{\textbf{Some definitions and main theorem}}\label{S3}

We consider binary Markov chains $\boldsymbol{\xi}= (\xi_{n})_{n=0}^{\infty}$ whose state space $X$ consists of two binary symbols, $X =\{0,1\}$, and for which Eq.~(\ref{e1}) holds. We assume that the chains $\bm{\xi}$ are time-homogeneous, which means that one-step transition probabilities $P(\xi_{n}=y|\xi_{n-1}=x)$ do not depend on time $n$, i.e., for some $\pi(x,y)$ Eq.~(\ref{e2}) holds; it is also assumed that some initial distribution of probabilities $P(\xi_{0} =x)$ on $X$ is given.

To proceed to formulation of our Theorem~\ref{T1}, we first present the notions of discrete capacity $\mathcal{C}$ and associated with this capacity thin and thick sets. The capacity $\mathcal{C}$ is assigned on $2^{\mathbb{N}}$; to define it, we consider binary codes of natural numbers. Let $k\in \mathbb{N}$, $(k\geq 1)$ and
$(\varepsilon_{0}, \ldots, \varepsilon_{p})$ be the binary code of $k$: $k = \sum_{i=0}^{p}\varepsilon_{i}2^{i}$ where $p\geq 0$, $\varepsilon_{i}\in \{0,1\}$ and $\varepsilon_{p} =1$ (binary expansion of $k$). Let $\nu(k)$ denotes the maximal of such $m$ ($0\leq m \leq p$), for which all the coefficients $\varepsilon_{i}$, $0\leq i\leq m$ of binary expansion of $k$ are equal to $1$.

\begin{definition}
For $e\subseteq \mathbb{N}$ we define
\begin{equation}\label{e3}
\mathcal{C}(e) = \sum_{k \in e}\nu(k).
\end{equation}
A set $e\subseteq \mathbb{N}$ is called thin (or, fine) set ($\mathfrak{F}$-set) if the relation
\begin{equation}\label{e4}
\sum_{p=1}^{\infty}2^{-p}\mathcal{C}(e\cap K_{p}) < \infty,
\end{equation}
where $K_{p}= \{k \in \mathbb{N}: 2^{p} \leq k < 2^{p+1}\}$, holds. If the set $e\subseteq \mathbb{N}$ is not thin (i.e.,  Eq.~(\ref{e4}) is failed), $e$ is called thick set ($\mathfrak{T}$-set).
\end{definition}

The $\mathcal{C}(e)$ from Eq.~(\ref{e3}) can be expressed in terms of binomial coefficients as follows. Let (for given $k\geq 1$) $\mu(k)$ denotes the maximal of such $m$ ($0\leq m \leq k$), for which all the binomial coefficients $\binom{k}{i}$, $0\leq i\leq m$ (first $m$ entries of $k$th line $(\binom{k}{0}, \binom{k}{1}, \ldots, \binom{k}{k})$ of the Pascale triangle), are odd numbers; one can prove that
$$
\mu(k) = 2^{\nu(k)}
$$
and, therefore,
$$
\mathcal{C}(e) = \sum_{k\in e}\log_{2}\mu(k).
$$

Since for infinite collection of bounded sets $e\subset \mathbb{N}$ and some positive constant we have $\mathcal{C}(e)\leq const.\mathcal{C}(\partial e)$ (cp.~\cite{4}; such inequality is mentioned also in \cite{9} when defining a capacity of clusters from $\mathbb{N}\times \mathbb{N}$, used in some models \cite{10}-\cite{11} of self-organized criticality), which is a characteristic property of classical capacities (e.g., \cite{12}), we call $\mathcal{C}$ a capacity. We note that $\mathcal{C}$ is differed from discrete capacity, considered in denumerable Markov chains and random walk (see, e.g., \cite{12}).

The next Proposition \ref{P1} contains some formal properties of capacity $\mathcal{C}$ and fine and thick sets (which we abbreviate  as $\mathfrak{F}$-sets and $\mathfrak{T}$-sets, respectively); we note that $\mathcal{C}(e) \geq 0$ for arbitrary $e\subseteq \mathbb{N}$.
\begin{proposition}\label{P1}
The next statements (a)-(f) are true: (a)~$\mathcal{C}(\emptyset) =0$ and $\mathcal{C}(\mathbb{N}) = \infty$.
(b)~If $e_{1}\subseteq e_{2}$ then $\mathcal{C}(e_1) \leq \mathcal{C}(e_2)$. (c)~ $\mathcal{C}(\{2^{p}\leq k <2^{p+1}\}) = (1+o(1))2^{p}$ \ $(p \to \infty)$. (d)~The $\mathbb{N}$ is  $\mathfrak{T}$-set. (e)~Every finite subset of $\mathbb{N}$ is $\mathfrak{F}$-set and finite union of $\mathfrak{F}$-sets is $\mathfrak{F}$-set.
(f) If $e$ is $\mathfrak{T}$-set and $e^{\prime}$ is $\mathfrak{F}$-set, then $e\cup e^{\prime}$ and $e\setminus e^{\prime}$ are  $\mathfrak{T}$-sets.
\end{proposition}

By using the next Proposition~\ref{P2} one can construct more complicated examples of thin and thick subsets of $\mathbb{N}$.

\begin{proposition}\label{P2}
Let for $p\geq 1$ the natural numbers $0\leq s_{p} \leq p$, $s_{p} \to \infty$  ($p\to \infty$) be given and $E \subseteq \mathbb{N}$ be defined as
\begin{equation}\label{s1}
E = \bigcup_{p=1}^{\infty}\{2^{p}\leq k < 2^{p+1}: \nu(k)\geq s_{p}\}.
\end{equation}
Then $E$ is $\mathfrak{F}$-set if and only if for $s_{p}$ the condition
$$
\sum_{p=1}^{\infty}\frac{s_{p}}{2^{s_{p}}} < \infty
$$
holds.
\end{proposition}

\begin{definition}
A number $a$ is called thick limit point ($\mathfrak{T}$-limit point or $\mathfrak{T}$-cluster point) of a given infinite numerical sequence $a_{k}$, $k\geq 0$ if there is a $\mathfrak{T}$-set $E\subseteq \mathbb{N}$  such that
$$
 \lim_{\substack{k\to \infty \\ k \in E}}a_{k} = a.
$$
A random binary sequence $\bm{\xi}= (\xi_{n})_{n=0}^{\infty}$ is called $\mathfrak{T}$-limit process for a given infinite series of random binary sequences $\bm{\xi}_{k} = (\xi_{n,k})_{n=0}^{\infty}$, $k\geq 0$ if for $x\in X$ and $n\geq 0$ the probability $P(\xi_{n}=x)$ is $\mathfrak{T}$-limit point for the sequence of probabilities $P(\xi_{n,k}=x)$, $k\geq 0$.
\end{definition}

The following Theorem~\ref{T1} is the main result of this paper.

\begin{theorem}\label{T1}
Let $\bm{\xi} = (\xi_{n})_{n=0}^{\infty}$  be time-homogeneous binary Markov chain for which Eq.~(\ref{E}) holds. Then the equi-distributed random binary sequence is the $\mathfrak{T}$-limit process for the sequence of higher-order differences $\bm{\xi}^{(k)} = (\xi_{n}^{(k)})_{n=0}^{^\infty}$, \ $k\geq 0$. More precisely, for $x\in X$ and $n \geq 0$ there is a $\mathfrak{T}$-set $E\subseteq  \mathbb{N}$ of the form (\ref{s1}) with $\sum_{p=1}^{\infty}s_{p}2^{-s_{p}} = \infty$,  for which
$$
\lim_{\substack {k\to \infty\\ k\in E}}P(\xi_{n}^{(k)} = x) =  \frac{1}{2}.
$$
\end{theorem}

In certain sense,  Theorem~\ref{T1} can be treated as the higher-order-difference version of the classical ergodic theorem for finite (two-state) Markov chains, where some notions from potential theory are now involved.

To the end, we present some characteristics of the sets $E$ from Theorem~\ref{T1} formulated in terms of their density in natural series. For $m\geq 1$ we denote $E_{m} = \{k\in E: 1\leq k\leq m\}$ and consider the ratio $\rho_{m}(E) = \dfrac{|E_{m}|}{m}$ where $|E_{m}|$ denotes the cardinality of $E_{m}$.

\begin{remark}\label{R1}
The sets $E\subseteq \mathbb{N}$ defined by Eq.~(\ref{s1}) in Proposition~\ref{P2} and presented in formulation of Theorem~\ref{T1} are of zero density in natural series: $\rho_{m}(E) \to 0$ as $m\to \infty$. The sets $E$ defined by Eq.~(\ref{s1}) can be such that the ratio $\rho_{m}(E)$ converges to $0$ as slow as we please: given $0< \delta_{m} \leq 1$, $\delta_{m} \downarrow 0$ the $\mathcal{T}$-set $E$ from Theorem~\ref{T1} can be constructed in such a way that $\rho_{m}(E) \geq \delta_{m}$ for all $m\geq 1$.
\end{remark}

\end{document}